\documentclass[preprint,12pt,fleqn]{elsarticle}

\usepackage{amssymb}
\usepackage{amsfonts}
\usepackage{amsthm}
\usepackage[english]{babel}
\usepackage{euscript}
\usepackage{bbm}
\usepackage{xfrac}
\usepackage{pb-diagram}

\usepackage{mathrsfs}
\usepackage{calrsfs}
\usepackage{amsbsy}
\newcommand{\fl}{\hspace*{-\mathindent}}
\newcommand{\textfrac}[2]{\textstyle{\frac{#1}{#2}}}

\newcommand{\D}{\EuScript{D}}
\newcommand{\E}{\mathrm{E}}
\newcommand{\J}{\mathrm{J}}

\usepackage{amssymb}
\usepackage{amsmath}
\usepackage{amsfonts}
\usepackage[english]{babel}
\usepackage{bbm}
\usepackage{xfrac}
\usepackage{color}

\usepackage{mathrsfs}
\usepackage{calrsfs}

\usepackage[all]{xy}
\SelectTips{cm}{}

\let\mathcal\mathscr

\let \kappa=\varkappa
\let \phi=\varphi

\begin{document}

\begin{frontmatter}

\title{Nonlocal conservation laws for the two-dimensional \\ Euler  equation in vorticity form}

\author{Oleg I. Morozov}   

\begin{abstract}
We combine the construction of the canonical conservation law and the nonlocal cosymmetry to derive a collection of  nonlocal conservation laws for the two-dimensional Euler equation in vorticity form. For computational convenience and simplicity of 
presentation of the results  we perform a complex rotation of the independent variables.
\end{abstract}


\end{frontmatter}


In this paper we consider the two-dimensional Euler equation in vorticity form \cite[\S~10]{LandauLifshits6}
\begin{equation}
\Delta u_t = \J (u, \Delta u),
\label{Euler_eq}
\end{equation}
where
$\Delta u = u_{xx} + u_{yy}$ and 
$\J(a,b) = a_x b_y-a_y b_x$. 
This  equation serves as a fundamental model in hydrodynamics and has been thoroughly investigated; for comprehensive discussions, see \cite{ArnoldKhesin} and references therein. We show that, in addition to the local conservation laws, the Euler equation admits 
nonlocal conservation laws ({\sc ncl}s), understood in the framework of the geometric theory of differential equations
\cite{KrasilshchikVinogradov1989,VK1999,KrasilshchikVerbovetskyVitolo2017}. Specifically, {\sc ncl}s are conservation laws whose coefficients depend on the pseudopotentials (the fiber variables) of a differential covering of the equation under study. We employ the construction of the canonical conservation law, see \cite{KrasilshchikVerbovetskyVitolo2017}, the results of \cite{Morozov2024}, in which a family of differential coverings for the  Euler equation was derived, and the nonlocal cosymmetry of equation \eqref{Euler_eq}, cf.  \cite{KrasilshchikMorozov2025}. We follow the definitions and the techniques of 
\cite{KrasilshchikVinogradov1984,KrasilshchikVinogradov1989,VK1999,KrasilshchikVerbovetskyVitolo2017}. All the considerations in the paper are local.

For a differential equation $\EuScript{E}$ defined as  $\{F=0\}$, the canonical conservation law   
is the cohomology class  of horizontal differential $(n-1)$-forms with a representative  $\Omega$ such that  
\begin{equation}
d_h\Omega = \left(\ell_F(q)\,p -q\,\ell^{*}_F(p)\right)\,dx^1 \wedge \dots \wedge dx^n,
\label{ncl_definition} 
\end{equation} 
where $\ell_F$  is the universal linearization operator associated with the function $F$, $\ell^{*}_F$ is the operator adjoint to $\ell_F$,
$d_h$ is the horizontal differential, and  $x^1$, ... , $x^n$ are the independent variables, see 
\cite[\S~6.1]{KrasilshchikVerbovetskyVitolo2017}. 
The form $\Omega$ is closed  on the Whitney sum of the  tangent and cotangent coverings of equation 
$\EuScript{E}$. This sum  is defined by the system
\[
\left\{
\begin{array}{rcl}
F &=& 0,
\\
\ell_F (q)&=& 0,
\\
\ell^{*}_F(p) &=& 0.
\end{array}
\right.
\]

\vskip 7 pt
\noindent
{\sc remark 1.}
In some cases, it is convenient to treat the variables $q$ and $p$ as odd (anticommuting), see \cite{KrasilshchikVerbovetskyVitolo2017} and references therein. Taking this into account, in \eqref{ncl_definition} and all formulas below, all instances of $q$ and its derivatives are placed to the left of $p$ and its derivatives.
\hfill $\diamond$ 
\vskip 7 pt

To simplify computations, we perform the following change of variables: we write equation
\eqref{Euler_eq}  as   
$\tilde{u}_{\tilde{t}\tilde{x}\tilde{x}}+\tilde{u}_{\tilde{t}\tilde{y}\tilde{y}} 
= 
\tilde{u}_{\tilde{x}}\,(\tilde{u}_{\tilde{x}\tilde{x}\tilde{y}}+\tilde{u}_{\tilde{y}\tilde{y}\tilde{y}})
-\tilde{u}_{\tilde{y}}\,(\tilde{u}_{\tilde{x}\tilde{x}\tilde{x}}+\tilde{u}_{\tilde{x}\tilde{y}\tilde{y}})$
and then put 
\begin{equation}
\tilde{t}=t, 
\quad
\tilde{x}=\textfrac{1}{2}\,(1+\mathrm{i})\,(x+y), 
\quad
\tilde{y} =-\textfrac{1}{2}\,(1-\mathrm{i})\,(x-y),
\quad
\tilde{u}=u, 
\label{complex_rotation}
\end{equation}
with $\mathrm{i} =\sqrt{-1}$. 
This yields the equation
\begin{equation}
\D (u_t) = \J (u, \D(u)),
\label{main_eq}
\end{equation}
where $\D = D_x\circ D_y$.

For equation \eqref{main_eq}, we have $F=\D (u_t) - \J (u, \D(u))$; hence the tangent and cotangent coverings are  obtained by appending the equations
\begin{equation}
\ell_F(q) = \D (q_t) - \J (q, \D(u)) -\J(u, \D (q))=0
\label{tangent_covering}
\end{equation}
and 
\begin{equation}
\ell^{*}_F(p) = -\D (p_t -\J(u,p)) - \J (\D(u), p)=0
\label{cotangent_covering}
\end{equation}
to equation \eqref{main_eq}.

In addition to the tangent and cotangent coverings, equation \eqref{main_eq} admits a  family of differential  coverings \cite{Morozov2024}
\begin{equation}
\left\{
\begin{array}{rcl}
s_t &=& \J(u, s)+\varepsilon\,\E(u),
\\
\J(\D(u), s) &=& \lambda + \mu\,\D(u)-\varepsilon\,\E(\D(u)),
\end{array}
\right.
\label{s_covering}
\end{equation}
where $\E(u) = x\,u_x+y\,u_y-2\,u$ and $\varepsilon, \lambda, \mu \in \mathbb{R}$.

Straightforward computations prove the following three assertions.

\vskip 7 pt
\noindent
{\sc Proposition 1.} 
{\it   The horizontal differential 2-form
\begin{equation}
\Omega_{q,p} = \D(q)\,p \,dx\wedge dy +K_1\,dy\wedge dt +K_2\,dt \wedge dx,
\label{CCL}
\end{equation}
where 
\[\fl
K_1 = 
u_y\,(q\,\D(p)+\D(q)\,p)
-\D(u_y)\,q\,p
-\textfrac{1}{2}\,(q_t\,p)_y
-q_y \,(p_t-u_x\,p_y)
+u_{yy}\,q\,p_x
\]
and 
\[\fl
K_2 = 
-u_x\,(q\,\D(p)+\D(q)\,p)
+\D(u_x)\,q\,p
+\textfrac{1}{2}\,(q_t\,p)_x
+q\,(p_{tx}-u_{xx}\,p_y)
-u_y\,q_x \,p_x,
\] 
defines the canonical conservation law of equation \eqref{main_eq}.
}
\hfill $\Box$

\vskip 7 pt

\noindent
{\sc Proposition 2.} 
{\it 
Let $s$ be a solution to equations \eqref{s_covering}. Then the function $p_0 = s-\lambda\,t\,x\,y-(\mu+2\,\varepsilon)\,t\,u$ is 
a solution to equation \eqref{cotangent_covering}.
}
\hfill $\Box$

\vskip 7 pt
\noindent
{\sc Proposition 3.}
{\it 
The Lie algebra of infinitesimal contact symmetries of equation \eqref{main_eq}  is generated by the functions 
\[
\begin{array}{rcl}
\varphi_1 &=& -t\,u_t-u,
\\
\varphi_2 &=& -u_t,
\\
\varphi_3 &=& -x\,u_x+y\,u_y,
\\
\varphi_4 &=& -t\,x\,u_x+t\,y\,u+x\,y,
\\
\varphi_5 &=& -x\,u_x-y\,u_y+2\,u,
\\
\varphi_6 (A_1)&=& -A_1(t)\,u_x+A_1^{\prime}(t)\,y,
\\
\varphi_7 (A_2)&=& -A_2(t)\,u_y-A_2^{\prime}(t)\,x,
\\
\varphi_8 (A_3)&=& A_3(t),
\end{array}
\]
where $A_i(t)$ are arbitrary smooth functions. 
}
\hfill $\Box$

\vskip 7 pt
Substituting  the cosymmetry $p=p_0$ from Proposition 2 and the generators $q=\varphi_i$ from Proposition 3  into the canonical conservation law \eqref{CCL}, we obtain a family of {\sc ncl}s $\Omega_{\varphi_i,p_0}$  for equation \eqref{main_eq}.   
The {\sc ncl}s corresponding to the generators $\varphi_j(A)$  with $j \in \{6,7,8\}$ depend on an arbitrary function  $A(t)$. 
The coefficients of $\Omega_{\varphi_i,p_0}$ depend on $t$, $x$, $y$, $u$, the derivatives of $u$ as well as  $s$ and derivatives of $s$ with respect to $t$, $x$, and $y$. We solve system \eqref{s_covering} with respect to $s_y$ and $s_t$, substitute the resulting expressions into $\Omega_{\varphi_i, p_0}$, and denote the result by $\omega_i$. In this way, we construct a collection of {\sc ncl}s 
$\omega_1$, ... , $\omega_8$, whose coefficients depend on $t$, $x$, $y$, $u$, the derivatives of $u$ up to and including fourth order, as well as on $s$, $s_x$, and $s_{xx}$.  
Using Green's formula, we can write  
\[
d_h \omega_i = (P_i\,F+Q_i\,G_1+R_i\,G_2)\, dt \wedge dx \wedge dy +d_h \eta_i  
\]
for each $i \in \{1, \dots, 8\}$, where $G_1$ and $G_2$  are the differences between the left and right sides of the first and second equations of system \eqref{s_covering}, respectively, $\eta_i$ are certain horizontal 2-forms, and all the coefficients $P_i$, $Q_i$, $R_i$ for $i \in \{1, \dots , 8\}$ are nonzero,  except for $Q_8 =0$. The explicit expressions for  the forms $\omega_i$  are quite lengthy, so we will restrict ourselves to the following two examples.  

\vskip 7 pt

\noindent
{\sc Example 1.}  
For $q  = \varphi_2 =u-x\,u_x$  
we have
\[
\omega_2 = -x\,\D(u_x)\, s \,dx\wedge dy
+\left( L_1\, s_{xx} +L_2\,s_{x}+L_3\,s+L_4\right)\,dy\wedge dt
\]
\[
\qquad\qquad
+\left(L_5\,s_{xx}+L_6\,s_{x}+L_7\,s+L_8\right)\,dt\wedge dx,
\] 
where 
\[
\fl
L_1 = 
u_y\,(u-x\,u_x)\,\D(u_y)\,\D(u_x)^{-1},
\]
\[\fl
L_2 = u_y\,(x\,u_x - u)\,\D(u_y)\,\D(u_{xx})\,\D(u_x)^{-2}
- u_y\,(x\,u_x - u)\,\D(\D(u))\,\D(u_x)^{-1}
\]
\[
+ \textfrac{1}{2}\,(x\,u_{tx} - u_t)\,\D(u_y)\,\D(u_x)^{-1}
- (x\,u_x - u)\,u_{yy}
- u_y\,(x\,\D(u) - u_y),
\]
\[\fl
L_3 =\textfrac{3}{2}\,x\,\J(u,\D(u))-u\,\D(u_y)-\textfrac{1}{2}\,u_{ty},
\]
\[\fl
L_4 =\lambda\,L_{4,1}+\mu\,L_{4,2}+\varepsilon\,L_{4,3},
\]
\[\fl
L_{4,1} = u_y\,(x\,u_x - u)\,\D(u_{xx})\,\D(u_x)^{-2}
+ \textfrac{1}{2}\,(x\,u_{tx} - u_t)\,\D(u_x)^{-1}
- \textfrac{3}{2}\,t\,x^2\,y\,J(u, \D(u))
\]
\[
+ t\,x\,y\,u\,\D(u_y)
+ x^2\,(t\,u_x - y)\,\D(u)
+ t\,y\,(x\,u_x - u)\,u_{yy}
- u_y\,(t\,u - x\,y)
\]
\[
- \textfrac{1}{2}\,t\,x\,(x\,u_{tx}- y\,u_{ty}-u_t),
\]
\[\fl
L_{4,2}= 
u_y\,(x\,u_x - u)\,\D(u)\,\D(u_{xx})\,\D(u_x)^{-2}
+ \textfrac{1}{2}\,(x\,u_{tx} - u_t)\,\D(u)\,\D(u_x)^{-1}
\]
\[
- \textfrac{3}{2}\,t\,x\,u\,J(u, \D(u))
+ t\,u^2\,\D(u_y)
- (x\,u + t\,(x\,u_t + u_y\,(u - 2\,x\,u_x)))\,\D(u)
\]
\[
- \textfrac{1}{2}\,(t\,x\,u_y\,u_{tx}
-t\,u\,u_{ty}
- u_y\,(t\,(3\,u_t - 2\,u_x\,u_y) - 2\,(x\,u_x - 2\,u)))
\]
\[
+ t\,u_x\,(x\,u_x - u)\,u_{yy},
\]
\[\fl
L_{4,3} = 
u_y\,(x\,u_x - u)\,(2\,\D(u) - y\,\D(u_y))\,\D(u_{xx})\,\D(u_x)^{-2}
\]
\[
+ y\,u_y\,(x\,u_x - u)\,\D(\D(u))\,\D(u_x)^{-1}
- \textfrac{1}{2}\,y\,(x\,u_{tx} - u_t)\,\D(u_y)\,\D(u_x)^{-1}
\]
\[
+ (x\,u_{tx} - u_t)\,\D(u)\,\D(u_x)^{-1}
- 3\,t\,x\,u\,J(u, \D(u))
+ 2\,t\,u^2\,\D(u_y)
\]
\[
+ \left(x\,(x\,u_x + y\,u_y - 4\,u) - 2\,t\,(x\,u_t + u_y\,(u - 2\,x\,u_x))\right)\,\D(u)
\]
\[
- \textfrac{1}{2}\,x\,(x + 2\,t\,u_y)\,u_{tx}
+ t\,u\,u_{ty}
+ 2\,t\,u_x\,(x\,u_x - u)\,u_{yy}
+ \textfrac{1}{2}\,(x + 6\,t\,u_y)\,u_t
\]
\[
- u_y\,(2\,u_x\,(x + t\,u_y) + y\,u_y - 5\,u),
\]
\[\fl
L_5=u_y\,(x\,u_x-u),
\]  
\[\fl
L_6=(x\,u_x-u)\,\D(u)-\textfrac{1}{2}\,(x\,u_{tx}-u_t)+x\,u_y\,u_{xx}, 
\]
\[\fl
L_7 =u\,\D(u_x)-\textfrac{1}{2}\,x\,u_{txx},
\]
\[\fl
L_8 =\lambda\,L_{8,1}+\mu\,L_{8,2}+\varepsilon\,L_{8,3},
\]
\[\fl
L_{8,1} = 
 \textfrac{1}{2}\,t\,(x^2\,y\,u_{txx}
+x\,y\,u_{tx}
- y\,u_t - 2\,x\,u_x^2)
+u_x\,(t\,u + x\,y)
- t\,x\,y\,u\,\D(u_x)
\]
\[
- t\,x\,(x\,u_x + y\,u_y - u)\,u_{xx}
- y\,u,
\]
\[\fl
L_{8,2} = 
 \textfrac{1}{2}\,t\,(x\,u\,u_{txx}
+ (3\,x\,u_x - 2\,u)\,u_{tx}
- u_t\,u_x)
+t\,(u_x\,(u - x\,u_x)\,\D(u)
- u^2\,\D(u_x))
\]
\[
+ t\,u_y\,(u - 2\,x\,u_x)\,u_{xx}
+ u_x\,(x\,u_x - u),
\]
\[\fl
L_{8,3} =
t\,(x\,u\,u_{txx}
+ (3\,x\,u_x - 2\,u)\,u_{tx}
- u_t\,u_x)
- 2\,t\,u^2\,\D(u_x)
+ 3\,u_x\,(x\,u_x - u)
\]
\[
+ (u - x\,u_x)\,(y + 2\,t\,u_x)\,\D(u)
+ (x\,(u - x\,u_x) + 2\,t\,u_y\,(u - 2\,x\,u_x))\,u_{xx}.
\]
Furthermore, we obtain 
\[\fl
d_h \,\omega_2 =
 ((x\,s_x+s-2\,t\,(\mu+2\,\varepsilon)\,u-2\,\lambda\,t\,x\,y)\,F
\]
\[\qquad
-x\,\D(u_x)\,G_1+(u-x\,u_x)\,G_2)\, dt\wedge dx\wedge dy +d_h \eta_2.
\]
\hfill $\diamond$

\vskip 7 pt
\noindent
{\sc Example 2.}  
When $q=\varphi_6=A(t)$,   the coefficient of $dx \wedge dy$ in \eqref{CCL}  vanishes, and we obtain the two-component {\sc ncl}. Two-component conservation laws are of importance in the study of recursion operators and other integrability structures for  multidimensional partial differential equations, see
\cite{KrasilshchikVerbovetskyVitolo2017,Krasilshchik2022,KrasilshchikVerbovetsky2022}.  
We get
\[
\omega_6 = \left((M_1\,s_{xx}+M_2\,s_{x}+M_{3}\,s+M_4)\,dy
\right.
\]
\[
\qquad\qquad 
\left.
-(M_5\,s_{xx}+M_6\,s_{x}+M_7\,s+M_{8})\,dx\right)\wedge dt,
\]
 with
\[\fl
M_1=A\,u_y\,\D(u_y)\,\D(u_x)^{-1},
\]
\[\fl
M_2= 
A\,(u_{yy}+u_y\,\J(\D(u),\D(u_x))\,\D(u_x)^{-2}) -\textfrac{1}{2}\,A^{\prime}\,\D(u_y)\,\D(u_x)^{-1},
\]
\[
\fl
M_3 = -A\,\D(u_y),
\]
\[\fl
M_4 =\lambda\,M_{4,1}+\mu\,M_{4,2}+\varepsilon\,M_{4,3}, 
\]
\[\fl
M_{4,1} =
A\,(t\,x\,y\,\D(u_y)-u_y\,L(u_{xx})\,\D(u_x)^{-2} -t\,(y\,u_{yy}+u_y))
+\textfrac{1}{2}\,A^{\prime}\,(t\,x-\D(u_x)^{-1}), 
\]
\[\fl
M_{4,2} =
A\,(t\,u\,\D(u_y) -t\,u_x\,u_{yy}+u_y\,(1-t\,\D(u))-u_y\,\D(u)\,D(u_xx)\,\D(u_x)^{-2})
\]
\[
+\textfrac{1}{2}\,A^{\prime}\,(t\,u_y-\D(u)\,\D(u_x)^{-1}),
\]
\[\fl
M_{4,3} =A\,(2\,t\,(u\,\D(u_y) -u_y\,\D(u)-u_x\,u_{yy}) - y\,u_y\,\J(\D(u),\D(u_x))\,\D(u_x)^{-2}+u_y
\]
\[
 -2\,u_y\,\D(u)\,\D(u_xx)\,\D(u_x)^{-2})
-\textfrac{1}{2}\,A^{\prime}\,(y\,\D(u_y)-2\,\D(u))\,\D(u_x)^{-1}
\]
\[
-\textfrac{1}{2}\,A^{\prime}\,(x+2\,t\,u_y), 
\]
\[\fl
M_5 =-A\,u_y, 
\]
\[\fl
M_6 =\textfrac{1}{2}\,A^{\prime}-A\,\D(u), 
\]
\[
\fl
M_7 = A\,\D(u_x),
\]
\[\fl
M_8 =\lambda\,M_{8,1}+\mu\,M_{8,2}+\varepsilon\,M_{8,3}, 
\]
\[\fl
M_{8,1} = 
A\,(t\,(x\,u_{xx}+u_x-x\,y\,\D(u_x))-y) -\textfrac{1}{2}\,A^{\prime}\,t\,y,
\]
\[\fl
M_{8,2} =
A\,(t\,(u_x\,\D(u)+u_y\,u_{xx}-u\,\D(u_x)-u_{tx})-u_x) -\textfrac{1}{2}\,A^{\prime}\,t\,u_x,
\]
\[
\fl
M_{8,3} =
A\,((y+2\,t\,u_x)\,\D(u)+(x+2\,t\,u_y)\,u_{xx}-3\,u_x -2\,t\,(u\,\D(u_x)+u_{tx}))
-A^{\prime}\,t\,u_x.
\]
Furthermore, we have
\[
d_h \,\omega_6 =
A\,(G_2 - (\mu+2\varepsilon)\,t\,F)\,dt\wedge dx \wedge dy +d_h\,\eta_6.
\]
\hfill $\diamond$

\vskip 7 pt

To obtain {\sc ncl}s for the Euler equation \eqref{Euler_eq}, one can apply the inverse of substitution 
\eqref{complex_rotation}  to the 2-forms $\omega_1$, ... , $\omega_8$. However, the computations are simpler if we rewrite equations \eqref{tangent_covering} and \eqref{cotangent_covering} for equation \eqref{Euler_eq}, replacing 
$\D$ with $\Delta$,  and then find the canonical conservation law  $\tilde{\Omega}_{q,p}$  for \eqref{Euler_eq}. 
For convenience in constructing the two-component conservation laws, we can take
\[
\tilde{\Omega}_{q,p} = \Delta q\,p\,dx \wedge dy +N_1\,dy \wedge dt +N_2\,dt \wedge dx,
\]
where
\[\fl
N_1 = 
q\,p_{tx} - q_x\,p_t 
+ u_y\,(\Delta\,q\,p -q_x\,p_x + q\,p_{xx})
+ \Delta \,q_y\,p  
- u_x\,(q\,p_{xy} +q_{xy}\,p) 
\]
\[
- u_{xy}\,(q_x\,p-\,q\,p_x) + u_{yy}\,q\,p_y
\]
and 
\[\fl
N_2 = 
q\,p_{ty} - q_y\,p_t
+u_x\,(q_y\,p_y-q\,p_y+q_x\,p_x-q_{yy}\,p)
+ u_y\,(q\,p_{xy} -q_y\,p_x)
\]
\[
- u_{xx}\,q\,p_x - u_{yy}\,q_x\,p - u_{xy}\,q\,p_y .
\]
Equation \eqref{Euler_eq}  admits a nonlocal cosymmetry 
$p_0=s-(\mu+2\,\varepsilon)\,t\,u-\textfrac{1}{4}\,\lambda\,(x^2+y^2)$
where $s$ is solution to system \eqref{s_covering} with $\D$ replaced by $\Delta$.  Substituting this cosymmetry in place of 
$p$,  the symmetries of equation \eqref{Euler_eq}  in place of $q$, and restricting on the solutions of system \eqref{s_covering} with $\D$ replaced by $\Delta$,  we obtain a collection of {\sc ncl}s  for this equation. While the coefficients of the resulting 2-forms turn out to be independent on $s_{xx}$,  their explicit expressions are too lengthy, so we omit them.

\section*{Acknowledgments}

I would like to express my sincere gratitude  to I.S. Krasil${}^{\prime}$shchik  for insightful discussions 
and valuable comments.

Computations  were done using the {\sc Jets} software \cite{Jets}.

\bibliographystyle{amsplain}

\end{document}